\documentclass[10pt]{amsart}
\title[Piecewise smooth Fermi-Ulam Models]
{Dynamics of some piecewise smooth\\ Fermi-Ulam Models}

\author{Jacopo De Simoi}
\address{Jacopo De Simoi\\
  Dipartimento di Matematica\\
  II Universit\`{a} di Roma (Tor Vergata)\\
  Via della Ricerca Scientifica, 00133 Roma, Italy.}
\email{{\tt desimoi@mat.uniroma2.it}}
\author{Dmitry Dolgopyat}
\address{Dmitry Dolgopyat\\
  Department of Mathematics\\
  University of Maryland\\
  4417 Mathematics Bldg,  College Park,  MD 20742, USA}
\email{{\tt dmitry@math.umd.edu}}
\urladdr{http://www.math.umd.edu/\~\ \hskip-4pt dmitry}

\thanks{We thank Anthony Quas and Carlangelo Liverani for useful discussions. This work was partially supported by the European Advanced Grant Macroscopic Laws and Dynamical Systems (MALADY)
  (ERC AdG 246953) and by NSF. Both authors are pleased to thank the Fields Institute in Toronto, Canada, for the excellent hospitality and
  working conditions provided in spring semester 2011.}
\usepackage{enumerate}
\usepackage{graphicx}
\usepackage{psfrag,amssymb,amsmath,amsthm}
\usepackage{verbatim}

\input{mapping.p}

\newcommand{\naturals}{\mathbb{N}}
\newcommand{\integers}{\mathbb{Z}}
\newcommand{\reals}{\mathbb{R}}

\newcommand{\torus}{\mathbb{T}}
\newcommand{\annulus}{\mathbb{A}}
\newcommand{\mes}{\text{mes}}
\newcommand{\EXP}{\mathbb{E}}
\newcommand{\Prob}{\mathbb{P}}

\newcommand{\ffrac}[2]{{#1/#2}}
\newcommand{\Orb}{\text{Orb}}
\newcommand{\continuous}[1]{\mathcal{C}^{#1}}

\newcommand{\bB}{\mathbf{B}}
\newcommand{\bG}{\mathbf{G}}
\newcommand{\bp}{\mathbf{p}}
\newcommand{\bT}{\mathbf{T}}
\newcommand{\bomega}{\boldsymbol{\omega}}

\newcommand{\brI}{\bar \aact}

\newcommand{\bF}{\mathbf{F}}
\newcommand{\brm}{\bar m}
\newcommand{\brv}{\bar v}
\newcommand{\brY}{\bar Y}
\newcommand{\brdelta}{\bar\delta}

\newcommand{\brLambda}{\bar\Lambda}
\newcommand{\brtau}{\bar\tau}
\newcommand{\brsigma}{\bar\sigma}
\newcommand{\brOmega}{\bar \Omega}
\newcommand{\cE}{\mathcal{E}}
\newcommand{\hsigma}{\hat\sigma}
\newcommand{\hPhi}{\hat{\Phi}}
\newcommand{\hPsi}{\hat{\Psi}}
\newcommand{\tmu}{\widetilde{\mu}}
\newcommand{\tnu}{\widetilde{\nu}}
\newcommand{\tfrm}{\widetilde{\frm}}

\newcommand{\RmII}{I\!\!I}

\newcommand{\eps}{\varepsilon}

\newcommand{\deh}{\textup{d}}

\newcommand{\st}{\textup{ s.t. }}
\newcommand{\bigo}[1]{\mathcal{O}\left(#1\right)}
\newcommand{\bigos}[2]{\mathcal{O}_{#2}\left(#1\right)}

\newcommand{\const}{\textup{Const}}

\newcommand{\tr}{\textup{Tr}}
\newcommand{\DDelta}{\Delta_1}
\newcommand{\Alg}{\mathbb{F}}
\newcommand{\sgn}{\text{sgn}}

\newtheorem{thm}{Theorem}
\newtheorem{lem}{Lemma}[section]
\newtheorem{prp}[lem]{Proposition}

\theoremstyle{remark}
\newtheorem{rmk}[lem]{Remark}
\theoremstyle{definition}

\newtheorem{exm}[lem]{Example}

\begin{document}
\begin{abstract}
  We find a normal form which describes the high energy dynamics of a class of piecewise smooth Fermi-Ulam ping pong
  models; depending on the value of a single real parameter, the dynamics can be either hyperbolic or elliptic. In the
  first case we prove that the set of orbits undergoing Fermi acceleration has zero measure but full Hausdorff
  dimension. We also show that for almost every orbit the energy eventually falls below a fixed threshold.
  In the second case we prove that, generically, we have stable periodic orbits for arbitrarily high energies,
   and that the set of Fermi accelerating orbits may have infinite measure.
\end{abstract}
\maketitle

\section{History and introduction}
In this paper we study the dynamics of piecewise smooth Fermi-Ulam ping pongs; Fermi and Ulam introduced such systems as
a simple mechanical toy model to explain the occurrence of highly energetic particles coming from outer space and
detected on Earth (the so-called \emph{cosmic rays}, see \cite{F49,F54}). The model describes the motion of a ball
bouncing elastically between a wall that oscillates periodically and a fixed wall, both of them having infinite
mass. Fermi and Ulam performed numerical simulations for the model and consequently conjectured (see \cite{U}) the
existence of orbits undergoing what is now called Fermi acceleration, i.e. orbits whose energy grows to infinity with
time; we refer to such orbits as \emph{escaping orbits}. Several years later, KAM theory allowed to prove that the
conjecture is indeed false. Namely, provided that the wall motion is sufficiently smooth, there are no escaping orbits
because invariant curves prevent diffusion of orbits to high energy (see \cite{LL, P1, P2}).  It was not many years (see
\cite{Z}) before the existence of escaping orbits was proved in some examples of piecewise-smooth motions; it is worth
noting that these examples were essentially the same that Fermi and Ulam were forced to investigate in their numerical
simulations, due to the relatively limited computational power they could use\footnote{The processing power of the
  1940's state-of-the-art computers used by Fermi and Ulam is about ten thousand times inferior than that of a low-end
  2010 smartphone}.  In this paper we study a more general class of piecewise smooth motions and we investigate
existence and abundance of escaping orbits in this setting.

Our main result is that, for all possible wall motions having one discontinuity, there is a parameter $\Delta$ which allows
to describe the dynamics of the ping pong for large energies. Moreover, there exists a sharp transition so that for
$\Delta \in (0,4)$ the system looks regular for large energies while for $\Delta \not\in [0,4]$ the system is chaotic
for large energies (see Figure \ref{f_tfrmElliptic} in Section \ref{s_Results}).
Similar phenomena happen in a wide class of piecewise smooth mechanical systems which for large
energies can be viewed as small ({\bf non-smooth}) perturbations of integrable systems such as, for example, the impact
oscillator \cite{DF}. However, in order to demonstrate the methods and techniques in the simplest possible setting, we
restrict our attention to the classical Fermi-Ulam model.

\section{Results.}
\label{s_Results}
We consider the following one-dimensional system: a unit point mass moves horizontally between two infinite mass walls,
between collisions the motion is free, so that kinetic energy is conserved, collisions between the particle and the
walls are elastic.  The left wall moves periodically, while the right one is fixed. The distance between the two walls
at time $t$ is denoted by $\el(t)$ which we assume to be strictly positive, Lipshitz continuous and periodic of period
$1$. It is convenient to study the orbit only at the moments of collisions with the moving wall.
Let $t$ denote a time of a collision of the ball with the moving wall, since $\el$ is periodic we take
$t\in\torus=\reals/\integers$. Let $v\in\reals$ be the velocity of the ball immediately after the collision.  Introduce
the notation $\annulus=\torus\times\reals$; the collision space is then given by
\[
\cspace=\{(t,v)\in\annulus\st v > -\dot\el(t)\}.
\]
We can thus define the \emph{collision map} $\cm:\cspace\to\cspace$
\begin{equation}\label{e_mapDefinition}
  \cm(t_n,v_n) = (t_n+\del(t_n,v_n),v_n-2\dot\el(t_n+\del(t_n,v_n))) = (t_{n+1},v_{n+1})
\end{equation}
where, for large\footnote{\emph{large} here means that the ball bounces off the fixed wall before the next collision
  with the moving wall} $v$, the function $\del$ solves the functional equation
\begin{equation}\label{e_defineDel}
  \del(t,v)=\frac{\ell(t) + \ell(t+\del(t,v))}{v}.
\end{equation}
It is a simple computation to check that the map $\cm$ preserves the volume form $\bomega=(v+\dot\el(t))\deh t\wedge\deh
v$.  Throughout this work we assume $\dot\el$ to be piecewise smooth with a jump discontinuity at $t=0$ only.  Define
the singularity line $\sing\subset\cspace$ as $\sing=\{ t=0 \}$ and $\sreg\in\cspace$ as the infinite strip of width
$\bigo{v^{-1}}$ bounded by $\sing$ and $\cm\sing$; introduce also $\psreg=\cmi\sreg$.

As a first step to study the dynamics of the mapping $\cm$ we describe the first return map of $f$ to the region
$\sreg$, which will be denoted by $\frm:\sreg\to\sreg$.  Our main result is a normal form for $\frm$ for large values of
$v$: introduce the notation $\el_0=\el(0)$, $\dot\el^\pm=\dot\el(0^\pm)$ and similarly for all derivatives. Define
$\Delta=\angnorm\el_0(\dot\el^+-\dot\el^-)$ and $\DDelta=\frac{1}{2}\angnorm^2\el_0^3(\ddot\el^+-\ddot\el^-)$ where
$\angnorm$ is given by
\begin{align}
  \label{DefAlphaDelta}
  \angnorm&=\int_0^1 \el^{-2}(s)\deh s.
\end{align}
We introduce a useful shorthand notation. Let $\psi\in\continuous{s}(A\subset\annulus)$; then we use the notation
$\psi=\bigos{v^{-k}}{s}$ to indicate that $v^k\psi$ is bounded for sufficiently large $v$ and the same is true for all
derivatives up to order $s$ included. For our analysis it is important to ensure that all sub-leading terms vanish
sufficiently fast for $v\to\infty$ along with all partial derivatives up to the fifth order.
\begin{thm}\label{t_firstReturn}
  There exist smooth coordinates $(\tau,\aact)$ on $\sreg$, such that
  the first return map of $\cm$ on $\sreg$ is given by
  \[
  \frm(\tau,\aact) = \hat\frm(\tau,\aact) +\frm_1(\tau,\aact) + r(\tau,\aact)
  \]
  where $\hat\frm(\tau, \aact)=(\bar\tau, \bar\aact)$ with
  \begin{align*}
    \bar\tau&=\tau-\aact\mod 1,&
    \bar\aact&= \aact+\Delta(\bar\tau-1/2),
  \end{align*}
  $\frm_1$ is a correction of order $\bigo{\aact^{-1}}$ of the form
  \[
  \frm_1(\tau, \aact)= \aact^{-1}(0, \DDelta((\bar\tau - 1/2)^2-1/12))
  \]
  and $r=\bigos{\aact^{-2}}5$. Finally $\bomega=\deh\tau\wedge \deh\aact. $
\end{thm}
Consequently, up to higher order terms, $\frm$ coincides with $\hat\frm$, where $\hat\frm$ is $\integers^2$-periodic in
appropriate ``action-angle'' variables and moreover $\deh\hat\frm=A$ is constant.  Thus $\hat\frm$ covers a map
$\tfrm:\torus^2\to\torus^2$; the map $\tfrm$ is known in the literature as the ``sawtooth map'' or the ``piecewise
linear standard map'' and it has been the subject of a number of studies, 
see e.g. \cite{BP, Bu, CM1, CM2, CDMMP, PV, W}. Notice that we have
\[ \tr(A)=2-\Delta. \]
Accordingly, $\deh\tfrm$ is elliptic if $\Delta\in (0,4)$ and it is hyperbolic otherwise.
\begin{exm}
Consider the case where velocity is piecewise linear, that is
\begin{equation}\label{e_FUM}
  \el_{A,B}(t)=B+A\left((t\mod{1})-\ffrac{1}{2}\right)^2.
\end{equation} 
This is  one of the cases which have been numerically
investigated in \cite{U}. We can choose the length unit so that $B=1.$
In this case $\el$ is positive for all $t$ for $A>-4.$
Remarkably, we can obtain an explicit expression for $\Delta(A)$, that is, for $\angnorm(A)$. Namely
\[
\angnorm(A)=\int_{0}^1(1+A(s-1/2)^2)^{-2}\deh s = 2|A|^{-1/2}\int_0^{|A|^{1/2}/2}(1+\sgn A\cdot\sigma^2)^{-2}\deh\sigma
\]
where $\sigma=|A|^{1/2}(s-1/2)$. Performing the integration we obtain
\[
\angnorm(A)=\frac{2}{A+4}+
\begin{cases}
  (|A|^{-1/2}/2)\log\frac{2+|A|^{1/2}}{2-|A|^{1/2}} & \text{if } -4<A\leq 0\\
  |A|^{-1/2}\arctan(|A|^{1/2}/2)&\text{if } A>0.
\end{cases}
\]
Recall that, by definition, $\Delta(A)=-2A(1+A/4)\angnorm(A)$. In particular, we find that $\Delta(-4)=4$ and
\begin{align*}
 \lim_{A\to-4^{+}}\Delta'(A) =\lim_{A\to-4^{+}} \frac{1}{2}\log\frac{2+|A|^{1/2}}{2-|A|^{1/2}}=+\infty.
\end{align*}
The graph of dependence of $\Delta$ on $A$ is shown on Figure \ref{f_Delta}. It shows that the dynamics
is hyperbolic for $A\in (-4, -a)$ where $a\approx -2.77927$ and $A>0$ and it is elliptic for remaining parameter values.
\begin{figure}[!h]
  \includegraphics[width=8cm]{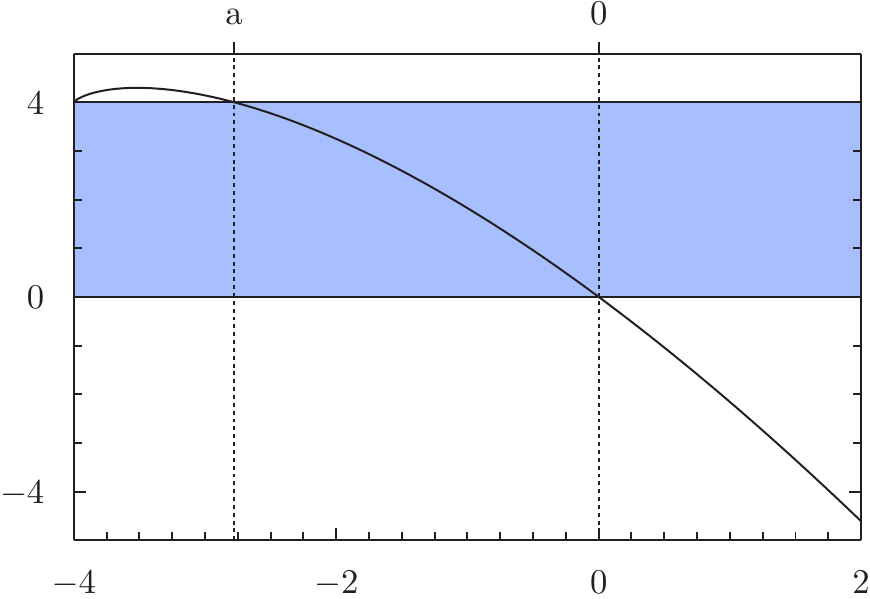}
  \caption{Graph of $\Delta$ as a function of $A$, where $\ell_{A,B}$ is given by \eqref{e_FUM} with $B=1$; the shaded area denotes the elliptic regime $\Delta\in(0,4)$.
  }
    \label{f_Delta}
\end{figure}
\end{exm}

\begin{figure}
  \includegraphics[width=9cm]{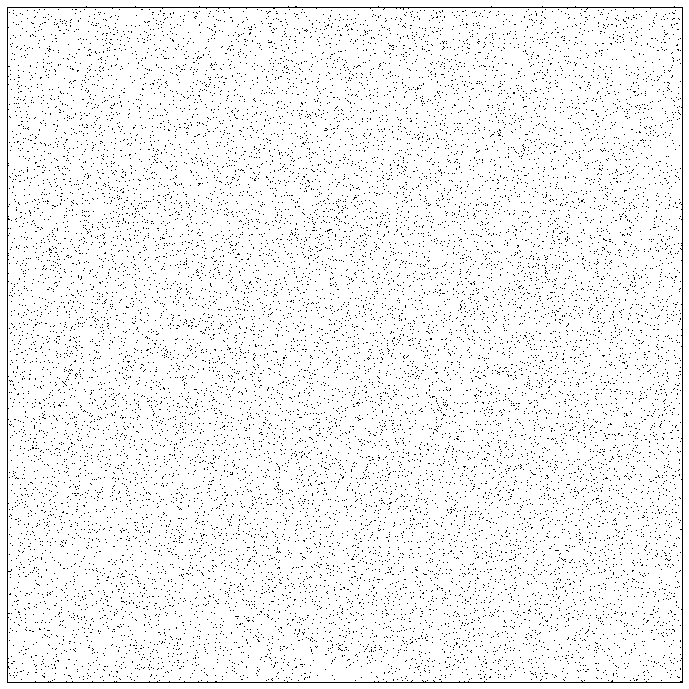}
  \vskip3mm
  \includegraphics[width=9cm]{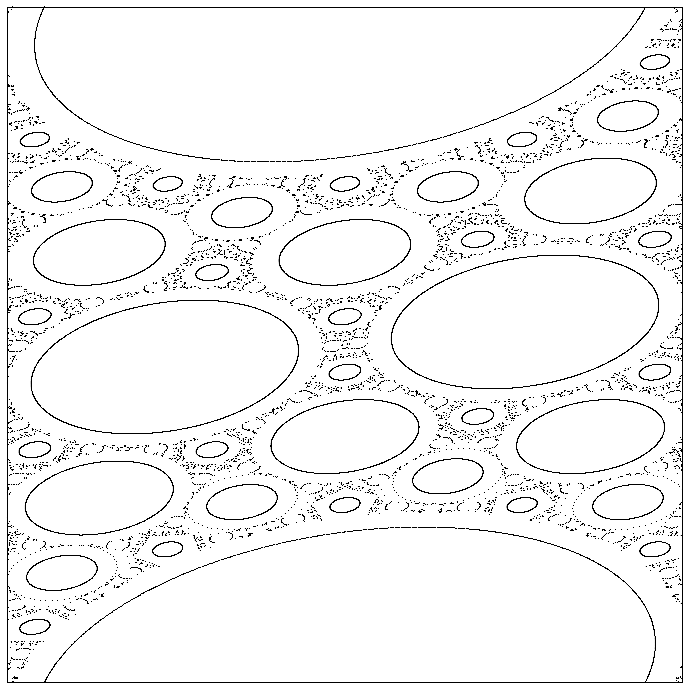}
  \caption{On the top: phase portrait of a single orbit of the map $\tfrm$ for $\Delta=-0.3$. On the bottom: phase
    portrait of selected orbits of the map $\tfrm$ for $\Delta=0.32$. Notice the strong prevalence of elliptic behavior;
    the ``chaotic'' region is given by forward and backward images of the singularity line.}
  \label{f_tfrmElliptic}
\end{figure}

Below we discuss the implications of the dichotomy between hyperbolic and elliptic regimes. Using the results of
\cite{C3} we obtain the following
\begin{thm}
  \label{t_ErgMix}
  If $|\tr(A)|>2$ then $\tfrm$ is ergodic, mixing and enjoys exponential decay of correlations for H\"older observables.
\end{thm}
On the other hand if $|\tr(A)|<2$ then $\tfrm$ is not ergodic. Namely, in this case, $\tfrm$ is a piecewise
isometry for the appropriate metric. Hence if $p$ is a periodic point of $\tfrm$, then a small ball around $p$ is
invariant by the dynamics. See Figure \ref{f_tfrmElliptic} for an example of phase portrait of $\tfrm$ in the two cases.

Note that if $p$ is periodic with period $N$ for $\tfrm$, it need not be periodic for $\hat\frm.$
In fact we have
\[ \hat\frm^N p=p+(0,n)\] for some $n\in\integers.$
If $n>0$ we say that $p$ is a \emph{stable accelerating orbit};
if $n<0$ we say that $p$ is a \emph{stable decelerating orbit}; finally 
if $n=0$ then $p$ is periodic for $\hat\frm.$

Consider for example the case $N=1$: we can find a periodic orbit $(\frac{1}{2}, 0)$;
if furthermore $\Delta>2$, we have a stable accelerating orbit
$(0, \frac{1}{2}+\frac{1}{\Delta})$ and stable decelerating orbit
$(0, \frac{1}{2}-\frac{1}{\Delta}).$ 

To analyze periodic points we can use the duality between the accelerating and decelerating periodic orbits. We have $\hat\frm=\tilt_\Delta\circ \frrm$ where
\begin{align}\label{e_frrm}
  \frrm(\tau, \aact)&=(\tau-\aact\mod 1, \aact),&
  \tilt_\Delta(\tau, \aact)&=(\tau, \aact+\Delta(\tau-\ffrac{1}{2})).
\end{align}
On the other hand $\hat\frm^{-1}(\brtau,\brI)=(\tau,\aact)$ with
\begin{align*}
  \aact&=\brI-\Delta(\brtau-\ffrac{1}{2}),&
  \tau&=\brtau+\aact.
  \intertext{Introducing $\sigma=1-\tau$ we rewrite the last equation as}
  \aact&=\brI+\Delta(\sigma-\ffrac{1}{2}),&
  \sigma&=\brsigma-\aact.
\end{align*}
In other words if $J$ denotes the involution $J(\tau, \aact)=(1-\tau, \aact)$ then
\begin{align*}
  (\tilt_\Delta\circ \frrm)^{-1}&=J \circ \frrm \circ \tilt_\Delta \circ J=\\
  &= J^{-1} \circ \frrm \circ \tilt_\Delta \circ J 
  = (\tilt_\Delta \circ J)^{-1} \circ (\tilt_\Delta \circ \frrm) \circ (\tilt_\Delta \circ J). 
\end{align*}
The existence of periodic orbits for other small periods is summarized in
table \ref{t_perAcc} (here we use parameter $\theta$ such that
$\tr(\deh\hat\frm)=2 \cos\theta, $ that is, $\deh\hat\frm$ is conjugated to a rotation
by $\theta$).

\begin{table}[!ht]
  \begin{tabular}{ccc}
    $N$ & periodic & accelerating/decelerating \\
    \hline
    1 & $(0, \pi)$ & $(\pi/2, \pi)$ \\
    2 & $(0, \pi/2)$ & - \\
    3 & $(0, \pi/3)$ & $(\pi/3,\pi/2)$\\
    4 & $(0, \pi/4)$ & $(3/4\pi,\pi)$\\
    5 & $(0, \pi/5)$ & $(\pi/2,3/5\pi)$\\
    6 & $(0, \pi/6) \cup (\pi/2,\pi)$ & $(5/6\pi,\pi)$\\
    7 & $(0, \pi/7)$ & $(3/7\pi,\pi/2)$\\
    8 & $(0, \pi/8)$ & $(7/8\pi,\pi)$\\
    \hline
  \end{tabular}
  \caption{Numerically observed ranges of parameters for which there exist a periodic or accelerating/decelerating orbit of ``period'' $N$; for $N\geq 9$ we could find intervals bounded by irrational multiples of $\pi$}\label{t_perAcc}
\end{table}
\begin{rmk}
  We believe that stable escaping orbits should exist for arbitrarily small positive values of $\Delta$, i.e. for each $\Delta_0>0$ there
  exists a $0<\Delta<\Delta_0$ such that the map $\hat\frm$ admits stable escaping orbits. However their ``period'' will necessarily
  to grow to infinity as $\Delta\to 0^+$; the smallest value of $\Delta$ for which we were able to find a stable
  escaping orbit is $\Delta=0.0916346$, for which we numerically obtained a period $501$ stable escaping orbit. 
\end{rmk}
We return now to the problem of energy growth; introduce the set of \emph{escaping orbits}
\[
\cE=\{(t_0, v_0): v_n\to\infty\}.
\]
\begin{thm}
  \label{t_Ell} If $|\tr(A)|<2$ then
  \begin{enumerate}[(a)]
  \item there exists a constant $C$ so that for each $\brv$ sufficiently large
    there exists an initial condition $(t_0,v_0)$ such that
    \begin{equation}\label{e_stabilitybrv}
      C^{-1}\brv< v_n <C\brv \text{ for all } n\in \integers.
    \end{equation}
    If additionally  $\DDelta\not=0$, then the same result holds for an adequately small ball around the point $(t_0,v_0)$.
  \item If $\tfrm$ has a stable accelerating
    orbit then $\mes(\cE)=\infty.$
  \end{enumerate}
\end{thm}
In particular $\mes(\cE)=\infty$ if $\Delta\in (\sqrt{3}, 4).$
\begin{thm}
  \label{t_EscHyp}
  If $|\tr(A)|>2$ then 
  \begin{enumerate}[(a)]
  \item $\mes(\cE)=0$.
  \item there exists a constant $C$ such that almost every orbit enters the region $v<C.$
    Moreover denote by $\bT$ the first time velocity falls below $C.$
    If we fix the initial velocity $v_0\gg 1$ and let the initial phase be random then 
    $\frac{\bT}{v_0^2}$ converges to a stable random variable of index $1/2,$ that is, there exists
    a constant $\bar{D}$ such that
    $$ P(\bT>\bar{D} v_0^2 t)\to \int_t^\infty \frac{e^{-1/2x}}{\sqrt{2\pi x^3}} \deh x \text{ as } v_0\to \infty.$$
  \end{enumerate}
\end{thm}
The proof of second part of the last theorem relies on the following result which is of independent interest.

\begin{thm}
  \label{t_CLT}
  Fix the initial velocity $v_0\gg 1$ and let the initial phase be random then 
  Fix $0<a<1<b.$ Consider the process defined by 
  $\bB^{v_0}(t)=\frac{v(v_0^2 t)}{v_0}$ 
  if $v_0^2 t$ is an integer with linear interpolation in between
  which is 
  stopped when velocity goes above $b v_0^2$ or below $a v_0^2.$
  Then, as  $v_0\to\infty$, $\bB^{v_0}(t)$ converges to a Brownian Motion started from 1
  and stopped when it reaches either $a$ or $b.$ 
\end{thm}

\begin{rmk}
  Note that $\bB^{v_0}(t)$ is equal to 
  $\frac{v(v_0^2 t)}{v_0}$ only if $v_0^2 t$ is an integer.
  It seems more natural to use this formula for all values of $t$ however this would lead to a different limit
  since as we shall see in the next section the ratio $v(n+\frac{1}{2})/v(n)$ has oscillations of order 1 while 
  $v(n+1)/v(n)-1$ is of order $1/v(n).$ 
\end{rmk}

Theorem \ref{t_CLT} makes Theorem \ref{t_EscHyp} plausible since the time  the Brownian Motion 
reaches a certain level has a stable distribution of index $1/2.$ However some work is needed to 
deduce Theorem \ref{t_EscHyp} from Theorem \ref{t_CLT} since the proof of Theorem \ref{t_CLT}
relies on a perturbative argument near $v=\infty$ whereas Theorem \ref{t_EscHyp} requires to 
handle small velocities as well since $v(\bT)\leq C.$ 

Theorem \ref{t_EscHyp} shows that the set of escaping orbits has zero measure so it is natural to ask about its Hausdorff dimension.  The next result extends the work \cite{dS} where a similar statement is proven for a smooth model of Fermi
acceleration. 

\begin{thm}
  \label{t_HD}
  If $\tr(A)>2$ then $\text{HD}(\cE)=2.$
\end{thm}
In other words, even though the set of escaping points is small from the measure theoretical point
of view, it is large from the point of view of dimension.
\section{The first return map}
If $\el$ were a smooth function, KAM theory would allow to conjugate the dynamics of $f$ for most initial conditions
with large values of $v$ to the dynamics of the completely integrable map
$\refm:\torus\times\reals^+\to\torus\times\reals^+$
\[
\refm:(\ang,\act)\mapsto(\ang+\act^{-1},\act).
\]
Consider the vertical line $S'\subset\torus\times\reals^+$ given by $S'=\{\ang=0\}$; let moreover $R'$ be the infinite strip of width $\bigo{\act^{-1}}$ bounded between $S'$ and $\refm S'$ i.e.
\[
R'=\{0\leq\ang<\act^{-1}\}.
\]
As a preparatory step we study the first return map of $\refm$ to the region $R'$.
\begin{prp}
  Let $\tau=\act\ang$ and consider coordinates $(\tau,\act)$ on $R'$. 
  Then the first return map of $\refm$ to the region $R'$ is given by the map $\frrm$ defined in \eqref{e_frrm}.
\end{prp}
\begin{proof}
  Let $k=\lfloor \act \rfloor$ and $\act = k+\hat\act$. We claim that
  \[
  \frrm(\tau,\act)=\begin{cases}
    \refm^k(\tau,\act)&\text{if }\tau\leq 1 - \hat\act\\
    \refm^{k+1}(\tau,\act)&\text{otherwise.}
  \end{cases}
  \]
  In fact, we can check by simple inspection that, denoting $\refm^k(\ang,\act)=(\ang_k,\act)$ we have
  \begin{align*}
    \ang_k&=\ang+\frac{k}{k+\hat\act}=\ang+ 1 - \frac{\hat\act}{\act},& \ang_{k+1}&=\ang+ 1 + \frac{1-\hat\act}{\act}
  \end{align*}
which implies our claim.
\end{proof}
In our systems, $\el$ is only piecewise smooth, consequently we expect to be able to define 
action-angle coordinates outside $\psreg$ only. 
\begin{lem}[Approximate reference coordinates]\label{l_actionangle}
  There exists a smooth coordinates change $\angact:(t,v)\mapsto(\ang, \act)$ such that if $(t,v)\not\in\psreg$,
  $\angact$ conjugates the collision map $\cm$ to the reference map $\refm$ up to high order terms
  \begin{equation}\label{e_estimatesApproximateMap}
    \refm - \angact\circ \cm\circ \angacti=(r_\ang,r_\act)
  \end{equation}
with $r_\ang=\bigos{v^{-4}}{5},\ r_\act=\bigos{v^{-3}}{5}$.
\end{lem}
\begin{proof}
  Recall the definition of $\angnorm$ given by \eqref{DefAlphaDelta} and introduce the notation $(t',v')=f(t,v)$. 
  Define the two functions (see e.g. \cite{Z} for a motivation of the formula defining $\ang$)
  \begin{align}\label{e_actionAngleVariables}
      \ang(t)&=\angnorm^{-1}\int_0^t \el^{-2}(s)\deh s, &\aact(t,v)=\angnorm\left[\int_t^{t'}{\el^{-2}(s)}\deh s\right]^{-1}.
    \end{align}
    It is immediate to observe that $\ang(t')=\ang(t) + \aact^{-1}(t,v)$. Since the expression defining $\aact$ is implicit, 
    we find it convenient to use a suitable approximation in our computations. 
    Define $\act:\annulus\setminus S\to\reals$ as 
    \begin{equation}\label{e_defAdia}
      2\angnorm^{-1} \act(\cdot,v)=v\el+\el\dot\el+\frac{1}{3}\frac{\el^2\ddot\el}{v}.
    \end{equation}
    We claim that $h:(t,v)\mapsto(\ang(t),\act(t,v))$ is the required change of coordinates. The first step is to obtain an approximate solution of \eqref{e_defineDel}.  Since $\el$ is Lipshitz continuous, we can find the solution by iteration. 
    Let $\delo0\equiv0$ and define for $n>0$
  \[
  \delo{n}(t,v) = \frac{\el(t)+\el(t+\delo {n-1}(t,v))}{v}.
  \]
  Then $\|\delo n-\delo{n-1}\|=\bigo{v^{-n}}$ and thus $\delo n\to\del$ uniformly.
  Consequently, if we express the solution as
  \begin{equation}\label{e_powerDel}
    \del = \sum_{n=1}^\infty \del_n \text{ with }\del_n(t,v)=\frac{a_n(t)}{v^n},
  \end{equation}
we can then find the functions $a_n$ by the previous argument.
  In particular, outside $\psreg$ we obtain that
\begin{align}\label{e_approxDelSmooth}
  \del_1(\cdot,v)&=\frac{2\el}{v},& 
  \del_2(\cdot,v)&=\frac{2\el}{v^2}\dot\el,&
  \del_3(\cdot,v)&=\frac{2\el}{v^3}(\dot\el^2+\el\ddot\el).
\end{align}
Assume now that $(t,v)\not\in\psreg$.
By expanding \eqref{e_actionAngleVariables} in Taylor series and using equations \eqref{e_approxDelSmooth} it is immediate to check that
\begin{equation}\label{e_errorIJ}
  \act=\aact+\bigos{v^{-2}}{5}.
\end{equation}
 Recall that 
  \begin{align*}
    r_\ang(t,v)&=\ang(t')-\ang(t)-\act(t,v)^{-1}, &    r_\act(t,v)&=\act(t',v')-\act(t,v).
  \end{align*}
  Thus estimate \eqref{e_errorIJ} immediately yields $r_\ang=\bigos{v^{-4}}{5}.$ The proof of 
  Lemma \ref{l_actionangle} is thus complete once we prove that 
  $$\act(t',v')-\act(t,v)=\bigos{v^{-3}}{5}. $$
  We begin by introducing a convenient notation. 
  Fix $(t,v)$. Recall that $r_\act=\act\circ f -\act$; denote $\act=\act(t,v),\ \act'=\act(t',v'),$ 
  $\el=\el(t),\ \el'=\el(t')$ and likewise for all derivatives.
  Notice that $\act v^{-1}$ is a polynomial in $v^{-1}$ with coefficients given by smooth functions of~$t$. Using \eqref{e_powerDel} we can express $\del$ in similar form, thus, by expanding in Taylor series the smooth function $\el$ and its derivatives we can write
  \[
  r_\act(t,v)=b_0(t)+\frac{b_1(t)}{v}+\frac{b_2(t)}{v^2}+r_\act^*(t,v)
  \]
  where $r_\act^*=\bigos{v^{-3}}{5}$. It amounts to a simple but tedious computation to show that our choice of $\act$ implies $b_0\equiv 0$, $b_1\equiv 0$ and $b_2\equiv 0$. Here we will only sketch the main steps of the computation.
   First, we obtain an expression for $\delta v$
  \[
  \delta v=v'-v =-2\el'= \delta v_0 + \delta v_1 + \delta v_2 + \bigo{v^{-3}} \text{ where}\]
  \begin{align*}
    \delta v_0 &= -2\dot\el &\delta v_1 &=- 4\frac{\el\ddot\el}{v} & \delta v_2 &= - 4\frac{\el\dot\el\ddot\el+\el^2\dddot\el}{v^2}.
  \end{align*}
Next, we expand $r_\act$ in Taylor series and collect terms of order $v^{-3}$ or higher in the function $r_\act^*$
  \begin{align*}
    r_\act&= \partial_t \act (\del_1+\del_2+\del_3) + \partial_v \act (\delta v_0+\delta v_1+\delta v_2) +\\
    &+ \frac{1}{2}\partial_{tt}\act\left(\del_1^2+2\del_1\del_2\right) + \partial_{tv}\act(\del_1\delta v_1+\del_1\delta v_2+\del_2\delta v_1)+\\
    &+ \frac{1}{6}\partial_{ttt}\act\del_1^3 + \frac{1}{2}\partial_{ttv}\act\del_1^2\delta v_1 + r_\act^*.
  \end{align*}
  Using the explicit form \eqref{e_defAdia} it is then simple to obtain
  \begin{align*}
    b_0&=\dot\el v\del_1+\el\delta v_0\\
    b_1&=\dot\el v^2\del_2+(\dot\el^2+\el\ddot\el)v\del_1+\frac{1}{2}\ddot\el v^2\del_1^2+\el v\delta v_1+\dot\el v\del_1\delta v_0.
\intertext{and finally}
    b_2&=\dot\el v^3\del_3+(\dot\el^2+\el\ddot\el)v^2\del_2+\frac{1}{3}(2\el\dot\el\ddot\el+\el^2\dddot\el)v\del_1+\\
    &+\ddot\el v^3\del_1\del_2 + \frac{1}{2}(3\dot\el\ddot\el+\el\dddot\el)v^2\del_1^2+\frac{1}{6}\dddot\el v^3\del_1^3+\\
    &+\el v^2\delta v_2 - \frac{1}{3}\el\ddot\el\delta v_0 + \dot\el v^2\del_1\delta v_1 + \dot\el v^2\del_2 \delta v_0 + \frac{1}{2}\ddot\el v^2\del_1^2\delta v_0.
  \end{align*}
  Now it is possible to conclude by substituting $b_j$ into the formulae obtained previously. 
\end{proof}
\begin{proof}[Proof of Theorem \ref{t_firstReturn}]
  It is simple to check by inspection that $(\tau,\aact)$ with $\tau=\aact\ang$ are smooth coordinates on $\sreg$ for sufficiently large $\aact$. 
  Let $(t,v)\in\sreg$, $(\bar t,\bar v)=\frm(t,v)$ and $(\tilde t,\tilde v)=\cmi(\bar t,\bar v)$. We use the convenient shorthand notation $\act=\act(t,v)$, $\tilde\act=\act(\tilde t,\tilde v)$ and $\bar\act=\act(\bar t,\bar v)$ 
  and similarly for $\ang$, $\tilde\ang$ and $\bar\ang$. By iteration of Lemma \ref{l_actionangle} we obtain
  \[
  \act(\tilde t,\tilde v) - \act(t,v)=\bigos{v^{-2}}5.
  \]
  We then claim that
  \begin{align}\label{e_tilt}
    \bar\act-\tilde\act&=\frac{1}{2}\angnorm(\dot\el^+-\dot\el^-)\left[\bar t \bar v (1-\frac{\smash{\dot\el^+}}{\bar v})-\el_0\right]+\\
    &+\frac{1}{4}\angnorm\frac{\ddot\el^+-\ddot\el^-}{\bar v}\left[(\bar t \bar v-\el_0)^2-\frac{1}{3}\el_0^2\right]
    +\bigos{\bar v^{-2}}5\notag.
  \end{align}
  In fact notice that
  \begin{align*}
    \ddot\el(\bar t)&=\ddot\el^++\bigo{\tilde v^{-1}}&
    \ddot\el(\tilde t)&=\ddot\el^-+\bigo{\tilde v^{-1}}\\
    \dot\el(\bar t)&=\dot\el^++\ddot\el^+\bar t+\bigo{\tilde v^{-2}}&
    \dot\el(\tilde t)&=\dot\el^-+\ddot\el^-\tilde t+\bigo{\tilde v^{-2}}\\
    \el(\bar t)&=\el_0+\dot\el^+\bar t+\frac{1}{2}\ddot\el^+\bar t^2 + \bigo{\tilde v^{-3}}&
    \el(\tilde t)&=\el_0+\dot\el^- \tilde t+\frac{1}{2}\ddot\el^-\tilde t^2 + \bigo{\tilde v^{-3}}
  \end{align*}
  and moreover
  \begin{align}\label{e_crucialtv}
    \tilde t &= \bar t - \frac{2\el_0}{\tilde v} + \frac{\dot\el^++\dot\el^-}{\tilde v}\bar t - \frac{2\el_0\dot\el^-}{\tilde v^2} + \bigo{v^{-3}}& \bar v = \tilde v - 2\dot\el(\bar t).
  \end{align}
By the definition of $\act$ we thus obtain
  \begin{align*}
    \bar\act -\tilde\act=\frac{1}{2}\angnorm&\Big[(\el(\bar t)-\el(\tilde t))\tilde v + \\
    &-(\el(\bar t)\dot\el(\bar t)) + \el(\tilde t)\dot\el(\tilde t)+\\
    &+\frac{1}{3}\frac{\el(\bar t)\ddot\el(\bar t)-\el(\tilde t)\ddot\el(\tilde t)}{\tilde v}\Big]+\bigos{\tilde v^{-2}}5
  \end{align*}
  from which \eqref{e_tilt} follows by a straightforward computation. Notice that by definition we have
  \begin{align*}
    2\angnorm^{-1}\bar \act&=\bar v\el_0+\dot\el^+\bar v\bar t+\el_0\dot\el^++\bigo{\bar v^{-1}}\\
    2\angnorm^{-1}\tilde \act&=\tilde v\el_0+\dot\el^-\tilde v\tilde t+\el_0\dot\el^-+\bigo{\tilde v^{-1}}.
    \intertext{Next, by definition of $\ang$}
    \bar t&=\angnorm\el_0^2\bar\ang(1 + \angnorm\el_0\dot\el^{+}\bar\ang)+\bigo{v^{-2}}
  \end{align*}
  from which we obtain
  \[
  \left(\bar t\bar v(1+\dot\el^+/{\bar v}) - \el_0\right) = 2\el_0(\bar\tau-1/2)+\bigo{v^{-2}}.
  \]
  Therefore we can rewrite $\bar \act$ as follows
  \[
  \bar \act = \act + \Delta(\bar\tau-1/2)+\frac{\DDelta}{\act}((\bar\tau-1/2)^2-1/12)+\bigos{\act^{-2}}5.
  \]
  Using estimate \eqref{e_errorIJ} we thus conclude that 
  $$\aact(\bar t,\bar v) - \aact(t,v)=\aact(\bar t,\bar v)+\bigos{\aact^{-2}}5. $$
  We now prove that 
  \begin{equation}
  \label{ChAng}
  \bar\tau = \tau - \act \mod 1.
  \end{equation}
   By Lemma \ref{l_actionangle} and definition we have
  \begin{align*}
    \tilde \act \tilde\ang &= \tau - \act -1& \bar \act \bar \ang &=\bar\tau.
  \end{align*}
  On the other hand, using the definition of $\ang$ and the approximate expressions for $\act$ given above, we obtain
  \begin{align*}
    \tilde \act\tilde\ang &=\frac{\tilde v+\dot\el^-}{2\el_0}\tilde t +\bigos{\act^{-2}}5& \bar \act\bar\ang &=\frac{\bar v+\dot\el^+}{2\el_0}\bar t+\bigos{\act^{-2}}5.
  \end{align*}
  From the last equation we obtain, using \eqref{e_crucialtv}, that
  \[
  \bar\act\bar\ang = \tilde \act \tilde \ang + 1 +\bigos{\act^{-2}}5.
  \]
Now  \eqref{ChAng} follows from \eqref{e_errorIJ}.
\end{proof}
\section{Hyperbolic case. Properties of the limiting map}
The proof of ergodicity of the map $\tfrm$ has been first given in \cite{V}.
Stronger statistical properties claimed by Theorem \ref{t_ErgMix}
  follow from the following general result.  Let $\bG$ be a piecewise linear hyperbolic automorphism of $\torus^2$ and denote by $S_+$ and $S_-$ the discontinuity curves of $\bG^{-1}$ and $\bG$, respectively; let $S=S_-\cup S_+$. 
  For any positive $n\in\naturals$ let $S_n=\bG^{n-1}S_+$ and $S_{-n}=\bG^{-(n-1)}S_-$; 
  assume for convenience $S_0=\emptyset$; let $S^{(n)}=\bigcup_{k=-n}^nS_k$.
\begin{prp}[Chernov, \cite{C3}]\label{t_chernov}
  Assume:
  \begin{enumerate}[(a)]
  \item $S_i\cap S_j$ is a finite set of isolated points if $i\not = j$;
  \item $S$ is everywhere transversal to the invariant stable and unstable directions;
  \item for every $m\geq 1$, the number of components of $S^{(n)}$ meeting at a single point is bounded by $K m$ for some constant $K$;
  \end{enumerate}
  then $\bG$ is ergodic, mixing and enjoys exponential decay of correlations for H\"older observables. 
\end{prp}
\begin{proof}[Proof of Theorem \ref{t_ErgMix}]
  If $\tr(A)>2$ then $\tfrm$ is a piecewise linear hyperbolic automorphism of $\torus^2$; we recall the explicit formula:
  \[
  \tfrm:(\tau, I)\mapsto(\tau-I\mod 1,I+\Delta((\tau-I\mod 1)-1/2)).
  \]
  Thus it is easy to check that $S_-$ is given by the diagonal circle $\tau=\aact$ and $S_+$ is given by the vertical
  circle $\tau=0$. It is then a simple linear algebra computation to prove that the stable and unstable slopes are given by the solution of the quadratic equation $h^2-\Delta h+\Delta=0$; thus we immediately obtain item (b) in the hypotheses of Proposition \ref{t_chernov}.
  Since $d\tfrm$ is constant at any point where it can be defined, 
  the $n$-th image of any line segment is a finite disjoint union of line segments parallel to each other; hence each point $p\in\torus^2$ can meet at most two of such segments, which proves item (c). Finally, unless the initial line segment is aligned to an invariant direction (stable or unstable), the slopes of line segments belonging to images at different times are different, which proves item (a) and concludes the proof.
\end{proof}

\section{Elliptic case. Growth of energy.}
\begin{proof}[Proof of Theorem \ref{t_Ell}.]
  In order to prove item (a) we will prove that for each $\brv$ sufficiently large there exists a stable periodic point
  $(t_*,v_*)\in\sreg$ whose orbit satisfies condition \eqref{e_stabilitybrv}; stability of the fixed point then implies
  that \eqref{e_stabilitybrv} holds for any initial condition $(t_0,v_0)$ in a small ball around $(t_*,v_*)$. We already
  noticed that the point $(\hat\tau,\hat\act)$ is a stable fixed point of the map $\hat\frm$ if $\hat\tau=1/2$ and
  $\hat\act\in\integers$; in order to prove existence of a stable fixed point of the first return map $\frm$ we would
  need to prove that the fixed point of $\hat\frm$ satisfies the non-degenerate twist condition. However, since
  $\hat\frm$ is piecewise linear, we actually need to consider the first return map $\frm$ as a $\bigo{\act^{-2}}$
  perturbation of the map $\smash{\bar\frm = \hat\frm+\frm_1}$ and check that $\bar\frm$ satisfies the twist
  condition. Since the perturbation term is small up to derivatives of sufficiently high order, we can conclude.
  
  Fix once and for all $\hat \act\in\naturals$ such that $|\act(0,\bar v)-\hat \act|<2$; let $\hat\lambda=\exp(i\hat\theta)$ be the multiplier at the fixed point $(\hat\tau,\hat\act)$ of $\hat\frm$. Since $\Delta\in (0,4)$ we have $\hat\lambda\not = 1$; then $\bar F$ will have a fixed point close to $(\hat\tau,\hat\act)$ that we denote by $(\bar\tau,\bar \act)$; introduce the parameter $\eps=\DDelta/I$; by inspection is is easy to see $\bar\act=\hat\act$ and $\bar\tau=\hat\tau+\eps/(12\Delta)+\bigo{\eps^{2}}$. Introduce coordinates $(\sigma,J)$ in a neighborhood of the fixed point $(\bar\tau,\bar\act)$ such that $(\tau,\act)=(\bar\tau+\sigma,\hat\act+\act)$. The expression for $\deh\bar\frm$ in these new coordinates is:
  \[
  \deh F(\sigma,J) = \left(
    \begin{array}{cc}
      1&-1\\
      \kappa+\eps\sigma&1-\kappa-\eps\sigma
    \end{array}
  \right)
  \]
  where $\kappa=\Delta+\eps^2/(6\Delta)$; denote by $\bar\lambda=\exp(i\bar\theta)$ the multiplier of the map $\bar\frm$ at $(\bar\tau,\bar\act)$: it is immediate to check that
  \begin{equation}\label{e_bartheta}
    \cos\bar\theta=\cos\hat\theta-\eps^2/(12/\Delta).
  \end{equation}
  In order to check the twist condition we perform a complex change of variables $(\sigma,J)\mapsto(z,\bar z)$ such that the map can be expressed as follows
  \[
  \bar F:z\mapsto z + \lambda z + A_3z^2 + A_4z\bar z + A_5\bar z^2.
  \]
  Then (see e.g. \cite{Churchill}) we need to ensure that:
  \[
  \Upsilon = 3|A_3|^2\frac{\bar\lambda+1}{\bar\lambda-1} + |A_5|^2\frac{\bar\lambda^3+1}{\bar\lambda^3-1} \not = 0;
  \]
  Notice that from the fact that $\bar F$ is symplectic we obtain that $|A_3|^2=|A_5|^2$; thus there are two possibilities for the twist condition to fail; either $A_3=A_5=0$ or $\bar\lambda$ solves the equation
  \[
  3\frac{\bar\lambda+1}{\bar\lambda-1} + \frac{\bar\lambda^3+1}{\bar\lambda^3-1}  = 0.
  \]
  It is easy to check that the above condition is given by either $\bar\theta=0$ or $\cos\bar\theta=-1/4$; both
  possibilities can be prevented by avoiding special values of $\eps$, according to \eqref{e_bartheta}, possibly by choosing
  a different $\hat\act$. Therefore, we just need to check that $A_3\not = 0$: from elementary linear algebra we 
  find that:
  \begin{align*}
    z &= \sigma + (1-\lambda)J&    z &= \sigma + (1-\bar\lambda)J;
  \end{align*}
  changing variables we obtain:
  \begin{align*}
    A_3&=-\frac{\lambda-1}{(\lambda-\bar\lambda)^2}\eps&
    A_4&=-\frac{1-\lambda-\bar\lambda}{\lambda-\bar\lambda}\eps&
    A_5&=\frac{\bar\lambda-1}{(\lambda-\bar\lambda)^2}\eps
  \end{align*}
  Which implies $|A_3|\not =0$ and concludes the proof of item (a).

  The proof of item (b) is analogous to the proof of the corresponding result obtained in \cite{DF}, Section 3, and will therefore be omitted.
\end{proof}
\section{Hyperbolic case. Measure of accelerating orbits.}
Note that part (a) of Theorem \ref{t_EscHyp} follows from part (b), however since the proof of part (b) is rather involved
we give a direct proof in this section. We expect Lemma \ref{l_ErgConserv} below to be useful for a wide range of mechanical systems. In particular, Theorem \ref{t_EscHyp}(a) is a direct consequence of 
Theorem \ref{t_ErgMix} and Lemma \ref{l_ErgConserv}.

Let $X$ be a Borel space and $Y$ be a subset of $X\times\naturals $ containing
$\{(x,m): m\geq \brm\}$ for some $\brm.$ Let $\Phi:Y\to Y$ be the map
$$ \Phi(x,m)=(\phi(x,m), m+\gamma(x,m)). $$
Assume that $\Phi$ is asymptotically periodic in the following sense. Denote
$T_k(x,m)=(x, m+k)$ and consider
$\Psi_k=T_k^{-1} \Phi T_k.$ Assume that there exist a map $\psi: X\to X$ preserving a probability measure
$\mu,$ and a function $\gamma:X\to \integers$ such that for each $M$ for each function $h$ supported
on $X\times [0, M]$ and each $l$ we have
\begin{equation}
\label{LimInfty}
 ||h\circ \Psi_k^l-h\circ \Psi^l ||_{L^2(\tmu)}\to 0 \text{ where }
\Psi(x,m)=(\psi(x), m+\gamma(x))
\end{equation}
and $\tmu$ is a product of $\mu$ and a counting measure on $\integers.$
Denote $(x_n, m_n)=\Phi^n(x_0, m_0)$ and let
$$ E=\{(x_0, m_0): m_n\to +\infty\}. $$ 
\begin{lem}
\label{l_ErgConserv}
Assume that
\begin{enumerate}[(i)]
\item $\psi$ is ergodic with respect to $\mu$;
\item $\int_X \gamma(x) \deh\mu(x)=0;$
\item $\Phi$ preserves a measure $\tnu$ with bounded density with respect to $\tmu;$
\item $||\gamma(x,m)||_{L^\infty}\leq K. $ 
\end{enumerate}
Then $\tnu(E)=0.$
\end{lem}
\begin{proof}%
By \cite{A} we know that conditions (i) and (ii) imply that $\Psi:X\times\integers \to X\times \integers$ is conservative.
That is for each subset $\brY$ of finite $\tmu$ measure the Poincare map $\hPsi:\brY\to\brY$ is defined almost
everywhere. Let $\brY=X\times [0, K+1]$ where $K$ is the constant from condition (iv). By Rohlin Lemma applied
to $\hPsi$ for each $\eps$ there exists a set $\Omega_\eps\subset \brY$ and a number $L_\eps$ such that
$$ \tmu(\Omega_\eps)< \eps \text{ and } \tmu((x,m): \hPsi^l(x,m)\not\in \Omega_\eps \text{ for }
0\leq l\leq L_\eps)<\eps. $$
In view of \eqref{LimInfty} and condition (iii) there exists a constant $C$ (independent of $\eps$) 
and a number $k(\eps)$ such that for $k\geq k(\eps)$ we have
$ \tnu(\brOmega_{k,\eps})<C\eps $ where
$$\brOmega_{k,\eps}=\{(x,m)\in T_k \brY: \Phi^l\not\in T_k \Omega_\eps \text{ for }l\in\naturals\}. $$
Let
$$ \Omega=\bigcup_{n\in\naturals} \left(T_{k(1/n^2)} \Omega_{1/n^2}\bigcup \brOmega_{k(1/n^2),(1/n^2)}\right). $$
Note that $\tnu(\Omega)<\infty.$ On the other hand if $(x,m)\in E$ then, due to condition (iv), its orbit $\Orb(x,m)$
visits $T_k\brY$ for all $k$ except for finitely many $k$. Hence $\Orb(x,m)\cap (\Omega\cap E)\neq\emptyset.$
Accordingly it suffices to show that $\tnu(\Omega\cap E)=0.$ However, by the foregoing discussion, the first return map
$\hPhi:(\Omega\cap E)\to(\Omega \cap E)$ is defined almost everywhere and by Poincare recurrence theorem for almost all
$(x,m)\in \Omega \cap E$ we have $\Orb(x,m)\cap X\times\{m\}\neq\emptyset.$ Thus for almost all points in $\Omega\cap E$
we have $(x,m)\not\in E.$ Therefore $\nu(\Omega\cap E)=0$ as claimed.
\end{proof}
\section{Hyperbolic case. Time of deceleration.}
\subsection{Plan of the proof.}
\label{SSPlan}
Here we prove Theorem \ref{t_EscHyp}(b). The argument of this section has 
many similarities with the arguments in \cite{CD2, D1, DSV} so we just 
indicate the key steps. 

The proof relies on the notion of {\it standard pair}.
A standard pair is a pair $\ell=(\gamma, \rho)$ there $\gamma$ is a curve such that $|\gamma|<1$
where $|\gamma|$ denotes the length of $\gamma,$ 
$\gamma'$ belongs to an unstable cone, 
$|\gamma''|\leq K_1,$ and $\rho$ is a probability density
on $\gamma$ satisfying $||\ln \rho||_{C^1(\gamma)}\leq K_2.$ We let $|\ell|$ denote the length of $\gamma.$
We denote by
$\EXP_\ell$ the expectation with respect to the standard pair 
$$\EXP_\ell(A)=\int_\gamma A(x) \rho(x) dx $$
and by $\Prob_\ell$ the associated probability measure, that is,
$ \Prob_\ell(\Omega)=\EXP_\ell(1_\Omega).$ 

An easy computation shows that if $I$ is sufficiently large on $\gamma$
then the standard pairs are invariant by dynamics, that is
$$ \EXP_\ell(A\circ\frm^n)=\sum_j c_j \EXP_{\ell_j}(A) $$
where $\sum_j c_j=1$ and $\ell_j$ are standard pairs. We need to know that most of $\gamma_j$ in this
decomposition are long. To this end let $r_n(x)$ be the distance from $x_n$ to the boundary of the component
$\gamma_j$ containing $x_n.$
\begin{lem}[Growth lemma]\label{l_Growth}$\phantom{1}$ %
\begin{enumerate}[(a)]
\item There exist constants $C>0$ and $\theta<1$ such that
$$\Prob_\ell(r_n(x)<\eps)\leq C\eps+\Prob_\ell(r_0(x)<\eps \theta^n).$$
\item There exists a constant $\eps_0$ such that if $n_0>K|\ln|\gamma||$ then
$$ \Prob_\ell(r_n(x)<\eps_0 \text{ for }  n=n_0,\ldots, n_0+k)\leq  C \theta^k. $$
\end{enumerate}
\end{lem}

The Growth Lemma is the key element of proving exponential mixing for $\tfrm$ (see \cite{C2, ChM})
and the argument used to prove the Growth Lemma for $\tfrm$ shows that this property is also valid for
small perturbations of $\tfrm.$

Given a point $x$ let $T_a$ be the first time $I_n<a$ if $I_0>a$ and be the first time 
$I_n>a$ if $I_0<a.$ Let $T_{a,b}$ be the first time then either
$I_n<a$ or $I_n>b.$ 

The proof of part (b) of Theorem \ref{t_EscHyp} depends on two propositions. The first one is an extended version
of Theorem \ref{t_CLT}. It will allow to handle large velocities. The second one gives an {\it a priori} bounded needed
to handle small velocities. 

Fix $0<a<1<b.$  Denote
\begin{equation}
\label{GK}
D^2=\sum_{n=-\infty}^\infty \iint_{\torus^2} A(x) A(\bF^n x) dx 
\end{equation}
where $\bF=G\circ T_\Delta$ and 
$G$ and $T_\Delta$ are defined by \eqref{e_frrm}.
\begin{prp}
\label{p_CLT}
Let $x$ be distributed according to a standard pair $\ell$ such that $I\sim I_0$ on $\ell$ and
$|\ell|>I_0^{-100}.$ Then
\begin{enumerate}[(a)]
\item The process $$B^{I_0}(t)=\frac{I_{\min(I_0^2 t, T_{a I_0, b I_0})}}{I_0}$$ 
converges to the Brownian Motion 
with zero mean and variance $D^2 t$
which is started from 1 and is stopped when it reaches either $a$ or $b;$
\item There exists $\delta>0$ such that $|\EXP_\ell(I_{T_{a I_0, b I_0}})-I_0|\leq C I_0^{1-\delta};$
\item There exists $\theta<1$ such that 
$\Prob_\ell(T_{a I_0, b I_0}>k I_0^2)\leq \max(\theta^k, I_0^{-100});$
\item Let $T^*_{a,b}=\min(T_{a I_0, b I_0}, I_0^3).$ Then 
  $$ \Prob_\ell(r_{T^*_{a,b}}(x)<\eps)\leq C I_0^3 \eps. $$
\end{enumerate}
\end{prp}

\begin{prp}
\label{p_Apr}
Given $\eps>0$ there exists $K(\eps)>0$ such that if $|\ell|>I_0^{-100}$ then
$$ \Prob_\ell(T_C> K(\eps) I_0^2)\leq \eps. $$
\end{prp}

Note that Proposition \ref{p_CLT} implies that 
$$ \Prob_\ell(T_{\delta I_0}\geq t I_0^2)\to \Prob(T^B_\delta\geq t)$$
where $T^B_\delta$ denotes the first time the Brownian Motion from Proposition \ref{p_CLT} reaches $\delta.$ 
Indeed 
$$ |\Prob_\ell(T_{\delta I_0}\geq t I_0^2)-
\Prob_\ell(T_{\delta I_0}\geq t I_0^2 \text{ and } T_{\delta I_0}\leq T_{A I_0})|\leq
\Prob_\ell(T_{\delta I_0} \geq T_{A I_0}). $$
By Proposition \ref{p_CLT} the RHS can be made as small as we wish by taking $A$ large.
On the other hand by Proposition \ref{p_CLT}
$$\Prob_\ell(T_{\delta I_0}\geq t I_0^2 \text{ and } T_{\delta I_0}\leq T_{A I_0})\to
\Prob_\ell(T^B_\delta \geq t \text{ and } T^B_\delta \leq T^B_A) $$
and the last expression can be made as close to 
$\Prob(T^B_\delta\geq t)$ as we wish by taking $A$ large. 

Next note that it is enough to prove Theorems \ref{t_EscHyp}(b) and \ref{t_CLT} with $v$ replaced by $I.$
Indeed, in view of \eqref{e_errorIJ}, \eqref{e_defAdia} and 
\eqref{DefAlphaDelta}, we have
\begin{equation}
\label{IV}
I\approx \frac{\angnorm \el(0)}{2} v
\end{equation}
which shows that $v$ can be replaced by $I$ in Theorem \ref{t_CLT}. Also \eqref{IV}
allows us to squeeze the first time $v$ goes below $C$ between the time $I$ goes
below $C_1$ and the time $I$ goes below $C_2$ and, in view of 
Proposition \ref{p_Apr},
the times to go below $C_1$ and $C_2$ satisfy the same estimates. 

We are now ready to derive part (b) of Theorem \ref{t_EscHyp} from 
Propositions \ref{p_CLT} and \ref{p_Apr}. 

We have
$$ \Prob_\ell(T_C\leq t I_0^2)\leq \Prob_\ell(T_{\delta I_0}\leq t I_0^2)\to 
\Prob(T^B_\delta\leq t)$$
and the last expression can be made as close to 
$\Prob(T^B_0 \leq t)$ as we wish by taking $\delta$ small.

Conversely
$$ \Prob_\ell(T_C\geq t I_0^2)\leq 
\Prob_\ell(T_{\delta I_0}\geq (t-K(\eps) \delta^2) I_0^2)+
\Prob_\ell(T_C(x_{T_{\delta I_0}})\geq K(\eps) \delta^2 I_0^2)$$
where $K(\eps)$ is given by Proposition \ref{p_Apr}.
The first term can be made as close to 
$\Prob(T^B_0\geq t)$ as we wish by taking small $\delta.$ To estimate the second term note that
$$ \Prob_\ell(T_C(x_{T_{\delta I_0}})\geq K(\eps) \delta I_0^2)$$
$$=\Prob(T_C(x_{T_{\delta I_0}})\geq K(\eps) \delta^2 I_0^2 \text{ and } r_{T_{\delta I_0}}(x)<(\delta I_0)^{-100})$$
$$+ 
\Prob(T_C(x_{T_{\delta I_0}})\geq K(\eps) \delta^2 I_0^2 \text{ and } r_{T_{\delta I_0}}(x)\geq (\delta I_0)^{-100})=I+\RmII.$$
Next
$$I\leq 
\Prob(r_{T_{\delta I_0}}(x)< (\delta I_0)^{-100}))=O(I_0^{-97})$$
by Proposition \ref{p_CLT}(d) and 
$$ \RmII=\Prob(r_{T_{\delta I_0}}(x)\geq I_0^{-100})) 
\Prob_\ell(T_C(x_{T_{\delta I_0}})\geq K(\eps) \delta^2 I_0^2|r_{T_{\delta I_0}}\geq (\delta I_0)^{-100}))$$
$$\leq
\Prob_\ell(T_C(x_{T_{\delta I_0}})\geq K(\eps) \delta^2 I_0^2|r_{T_{\delta I_0}}\geq (\delta I_0)^{-100}))\leq \eps
$$
where the last inequality follows by definition of $K(\eps).$

This completes the derivation of Theorem \ref{t_EscHyp}(b) from Propositions \ref{p_CLT} and \ref{p_Apr}.
It remains to establish the propositions. Proposition \ref{p_CLT} is proven in section
\ref{ss_CLT} and Proposition \ref{p_Apr} is proven in section \ref{ss_Apr}.

\subsection{Central Limit Theorem.}
\label{ss_CLT}
Let 
$$\frm_{I}^\dagger(\aact, \tau)=\hat\frm(\aact, \tau)+[\aact]^{-1}(0, \DDelta((\bar\tau - 1/2)^2-1/12)).
$$
Note that $\frm_I^\dagger$ approximates $\frm$ up to error $O(I^{-2}).$
Next consider a mapping of the $\torus^2$ given by
$$\brfrm_N(\aact, \tau)=\tfrm(\aact, \tau)+N^{-1}(0, \DDelta((\bar\tau - 1/2)^2-1/12)).$$
Then $\frm_I^\dagger$ 
locally covers $\tfrm_{[I]}$; also $\brfrm_N$ preserves the measure
$\deh\aact \deh\tau.$ The proof of Theorem \ref{t_ErgMix} shows that $\brfrm_N$ is exponentially mixing.
In particular, if $|\ell|\geq \eps_0$ then
$$ \EXP_\ell(A\circ \brfrm_N^n)=\iint_{\torus^2} A \deh\aact \deh\tau+O(\theta^n).$$
We use this property to establish the following estimate

\begin{lem}[Averaging Lemma]
\label{l_AV}
Suppose that $|\ell|>I_0^{-100}.$ Let $n=K\ln I_0$ where $K$ is sufficiently large.
Let $A$ be a piecewise smooth periodic function.
\begin{enumerate}[(a)]
\item $\EXP_\ell(A\circ \frm^n)=\iint_{\torus^2} A \deh\aact \deh\tau+O(I_0^{-2+\delta});$
\item There is $L>0$ such that
  $$\EXP_\ell(A(\frm^n x) A(\frm^{n+k} x))=\iint_{\torus^2} A(x) A(\tfrm^k x) \deh\aact \deh\tau+O(I_0^{-\beta} L^k).$$
\end{enumerate}
\end{lem}
The proof of this lemma is similar to the proof of Proposition 3.3 in \cite{CD1}. 
The proof of part (a) proceeds in two steps.
First, if $|\ell|>\eps_0$ then we use the shadowing argument to show that
\begin{equation}
\label{Shad}
\EXP_\ell(A\circ\frm^n)=\EXP_\ell(A\circ \brfrm_{[I_0]}^n)+O(I_0^{-2+\delta})
\end{equation}
and then use exponential mixing of $\brfrm_{[I_0]}.$ In the general case we find a function
$n(x)<\frac{K}{2} \ln I_0$ such that
$$ \EXP_\ell(A(\frm^{n(x)} x))=\sum_j c_j \EXP_{\ell_j}(A) $$ 
and $\sum_{|\ell_j|\leq \eps_0} c_j \leq I_0^{-100}$ and then 
apply \eqref{Shad} to all long components $\ell_j.$ 

To prove part (b) we first use the foregoing argument to show that
$$\EXP_\ell(A(\frm^n x) A(\frm^{n+k} x))=\iint_{\torus^2} A(x) A(\frm^k x) \deh\aact \deh\tau+O(I_0^{2-\delta} L^k)$$
(the factor $L^k$ accounts for the exponential growth of the Lipshitz norm of $A(A\circ \frm^k)$)
and then use the shadowing argument again to show that
$$ \iint_{\torus^2} A(x) A(\frm^k x) \deh\aact \deh\tau=\iint_{\torus^2} A(x) A(\tfrm^k x) \deh
\aact \deh\tau
+O(I_0^{-\beta}). $$
It is shown in \cite{CD2}, Appendix A that Lemmas \ref{l_Growth} and \ref{l_AV} imply parts (a) and (b)
of Proposition \ref{p_CLT}. 
We note that the error bound $$O(I^{-(2-\delta)})\ll I^{-1}$$ is needed to compute the drift of the
limiting process; to compute its variance it is enough that 
$\frm=\hat\frm+O(I^{-1})$ and that $\hat\frm$ covers $\tfrm$ which satisfies the CLT in the sense that
$$ \frac{I_n}{\sqrt{n}}\Rightarrow \text{Normal}(0, D^2)$$
where the diffusion coefficient $D^2$ is given by the Green-Kubo formula \eqref{GK}.
(In fact \eqref{GK} is the Green-Kubo formula for $\bF=G\tfrm G^{-1}$ but 
$\bF$ and $\tfrm$ clearly have the same transport coefficients.)

Next, part (a) of Proposition \ref{p_CLT} implies part (c) with $k=1,$ that is, there is $\theta<1$
such that
\begin{equation}
\label{CanEscape}
 \EXP_\ell(T_{aI_0, b I_0}<I_0^2)\leq \theta. 
\end{equation}
For $k>1$ we argue by induction applying \eqref{CanEscape} to all long components of $\frm^{(k-1) I_0^2} \gamma$
which have not escaped by the time $(k-1) I_0^2.$ 
Finally 
$$ \Prob_\ell(T^*_{a,b}\leq \eps)\leq \sum_{m=K\ln I_0}^{I_0^3} 
\Prob_\ell(r_m(x)<\eps) $$
so part (d) follows from part (a) of Lemma \ref{l_Growth}.
 
\subsection{A priori bounds for the return time}
\label{ss_Apr}
Let $\sigma_0$ be the first time when $|I_\sigma-2^{m_0}|\leq \Delta.$ For $j\geq 1$ we define
$\sigma_j$ inductively as follows. Assume that $\sigma_{j-1}$ was already defined so that 
$|I_{\sigma_{j-1}}-2^{m_{j-1}}|\leq \Delta.$ Let $\hsigma_j$ be the first time after $\sigma_{j-1}$ when
either $|I_\sigma-2^{m_j+1}|\leq \Delta$ or $|I_\sigma-2^{m_j-1}|\leq \Delta.$ Let
$\sigma_j=\min(\hsigma_j, \sigma_{j-1}+2^{3 m_{j-1}}).$ If either 
$\hsigma_j\geq \sigma_{j-1}+2^{3 m_{j-1}}$ or $r_{\sigma_j}<2^{-100 m_j}$ or
$2^{m_j}<\brI$ then we stop otherwise we continue and proceed to define $\sigma_{j+1}.$ If we stop we let
$j^*=j$ be the stopping moment. If we stop for the first or the second reason we say that we have
an emergency stop, otherwise we have a normal stop. By the discussion at the end of section \ref{SSPlan}
the lower cutoff in Proposition \ref{p_Apr} is not important so to prove the proposition it is enough to control
the first time when $I_n$ is close to 
$2^{\brm}$ with $2^{\brm}<\brI.$ In other words we need to control
$\sigma_{j^*},$ especially if it is a normal stop. Also since $\sigma_0$ is unlikely to be large by part (c)
of Proposition \ref{p_CLT} (in fact, part (a) would also be sufficient for our purposes) we need to control
$\sigma_{j^*}-\sigma_0.$ 

Let $\Alg_j$ be the $\sigma$-algebra generated by
$(m_0, \sigma_0),\cdots,(m_j, \sigma_j).$ Proposition \ref{p_CLT} implies that
\begin{align*}
\Prob_\ell(m_{j+1}=m_j+1|\Alg_j)&=\frac{1}{3}+o(1), &\brI&\to\infty\\
\Prob_\ell(m_{j+1}=m_j-1|\Alg_j)&=\frac{2}{3}+o(1), &\brI&\to\infty\\
\Prob_\ell(\sigma_{j+1} \text{ is an emergency stop}|\Alg_j)&=\bigo{2^{-97 m_j}}.
\end{align*}
Let
$\xi_j$ be a random walk with $\xi_0=m_0$ and 
\begin{align*}
\Prob(\xi_{j+1}=\xi_j+1)&=0.4, &\Prob(\xi_{j+1}=\xi_j-1)&=0.6
\end{align*}
Let $\Lambda_j$ be iid random variables independent of $\xi$s such that
$$ \Prob(\Lambda_j=k)=\begin{cases} K \theta^k & \text{ if }k\geq k_0 \\
                                                              0 & \text{otherwise} 
                                       \end{cases} $$
where $k_0$ is sufficiently large and $K=\frac{1-\theta}{\theta^{k_0}}.$
Let $\brLambda_j=\min(2^{2 \xi_j} \Lambda_j, 2^{3 \xi_j}).$
Proposition \ref{p_CLT} allows us to construct a coupling such that for $j\leq j^*$
\begin{align*}
  m_j&\leq \xi_j, &\sigma_j&\leq \sigma_0+\sum_{m=0}^{j-1} \brLambda_m.
\end{align*}
Now a standard computation with random walks shows that Proposition \ref{p_Apr} is valid for the 
random walk itself. Consequently, given $\eps$, there exists $K(\eps)$ such that
$$ \Prob_\ell\left(\sigma_{j^*}-\sigma_0\geq 2^{m_0} \frac{K(\eps)}{2}\right)\leq \frac{\eps}{10}. $$
Unfortunately $j^*$ need not to be a normal stop, it can be an emergency stop as well. 
To deal with this problem let 
$$\bp_k=\Prob_\ell(j^* \text{ is an emergency stop and }m_{j^*}=k).$$
Denote $\Omega_{kl}=\{j^*$ is an emergency stop, $m_{j^*}=k$ and $j^*$ is the $l$-th visit to $k\},$
$V_{kl}=\{k$ is visited at least $l$ times $\}.$ Then
$$ \bp_k\leq \sum_{l=1}^\infty \Prob_\ell (\Omega_{kl})=
\sum_{l=1}^\infty \Prob_\ell (V_{kl}) \Prob_\ell (\Omega_{kl}|V_{kl}). $$
By Proposition \ref{p_CLT}(d) 
$$\Prob_\ell(\Omega_{kl}|V_{kl})\leq C 2^{-97 k} $$
while the existence of the coupling with the random walk discussed above implies that
$$ \Prob_\ell(V_{kl} )\leq \theta^l.$$
Therefore 
\begin{equation}
\label{Stopk}
\bp_k\leq C 2^{-97 k}. 
\end{equation}
Accordingly, by choosing $\brI$ large enough we can make the probability of an emergency stop
less than $0.1.$ However we can not decrease that probability below $\eps/2$ if $\brI$ is fixed, so more work is needed.
 
First, we note that 
$\Prob_\ell(m_{j^*}>m_0/2)=O(I_0^{-50})$ so it can be neglected. 
Secondly, an argument similar to one leading to \eqref{Stopk} shows that 
$$ \Prob_\ell(r_{\sigma_{j^*}}< I_0^{-100})=O(I_0^{-50}). $$
Next, if $r_{\sigma_{j^*}}> I_0^{-100},$ $m_{j^*}<m_0/2$ and $j^*$ is an emergency stop let
$\brsigma$ be the first time after $\sigma_{j^*}$ such that $r_{\brsigma}(x)\geq \eps_0.$ 
By the Growth Lemma
$$ \Prob_\ell(\brsigma-\sigma_{j^*}>K \ln I_0)\leq I_0^{-100}$$
if $K$ is large enough.

If $\brsigma-\sigma_{j^*}<K \ln I_0$ then we can repeat the procedure described above with
$x$ replaced by $x_{\brsigma}.$ If the second stop is a normal one we are done, otherwise we try the third time
and so on. We have
$$ \Prob_\ell(\text{First }k \text{ stops are emergency stops})\leq (0.1)^k $$
which can be made less than $\eps/10$ if $k$ is large enough.
Next we have
$$ \Prob_\ell(k^*<k, T_{\brI} \leq K(\eps) I_0^2)\leq \frac{\eps}{5} $$
since the first try takes less than $\frac{K(\eps)}{2} I_0^2$ with probability greater than $1-\frac{\eps}{10}$
and all other tries take time $O(I_0)$ since with overwhelming probability we start those tries below level
$O(\sqrt{I_0}).$ This concludes the proof of Proposition \ref{p_Apr}.

\section{Hyperbolic case. Dimension of accelerating orbits.}
\begin{proof}[Proof of Theorem \ref{t_HD}] Foliate the phase space by line segments parallel 
to the unstable direction
of the limiting map $\hat\frm.$ It suffices to show that, given $s<1$ there exists $\brI$ such that if
$\Gamma$ is a leaf of our foliation and $I\geq \brI$ on $\gamma$ then
$\text{HD}(\Gamma\cap \cE)>s.$ 

By Theorem \ref{t_ErgMix} the limiting map satisfies CLT. That is, for any unstable curve $\gamma$, if the initial conditions
are distributed uniformly on $\gamma$ then 
$\frac{\hat I_n-I_0}{\sqrt n}$ converges to a normal distribution with zero mean and some variance $D$
(here we are using the notation $\hat\frm^n(\tau_0,I_0)=(\hat\tau_n,\hat I_n)$).
In particular there exists a constant $\kappa>0$ such that, for sufficiently large $n_0$, we have
\begin{equation}
\label{CLTUp}
\Prob_\ell(\hat I_{n_0}-I_0>\kappa\sqrt{n_0})>\frac{1}{3}
\end{equation}
where $\ell$ denotes the standard pair $(\gamma, \const).$
Moreover, given $\delta<\brdelta$, we can find $n_0$ so that \eqref{CLTUp} holds uniformly for all 
curves of length between $\delta$ and $\brdelta.$ 
Let $\hat r_n(x)$ denote the distance from $\hat\frm^n x$ to the boundary
of the component of $\hat\frm^n \Gamma$ containing  $\hat\frm^n x.$ 
By the Growth Lemma (Lemma \ref{l_Growth}) if $\delta$ is sufficiently small than for sufficiently large $n_0$ we have
$$ \Prob_\ell(\hat r_{n_0}(x)<3\delta)<\frac{1}{10} $$
provided that $\gamma$ is longer than $\delta.$
By Theorem \ref{t_firstReturn} we can take $\brI$ so large that if $I_0>\brI$ on $\gamma$ then 
$$ \Prob_\ell(I_{n_0}-I_0>\kappa\sqrt{n_0} \text{ and } r_{n_0}(x)>3\delta)\geq \frac{1}{5} $$
where $\frm(I_0, \tau_0)=(I_n, \tau_n)$ and $r_n(x)$, as before, denotes the distance from $\frm^n x$ to the boundary of
the component of $\frm^n \Gamma$ containing $\frm^n x.$ Note that any curve of length greater than $3\delta$ can be
decomposed as a disjoint union of curves with lengths between $\delta$ and $2\delta.$ Hence $\frm^{n_0} \gamma\supset
\bigcup_j \gamma_j $ where on each $\gamma_j$ the action grew up by at least $\kappa n_0$ and the total measure of
$\bigcup_j \frm^{-n_0} \gamma_j$ is at least $\mes(\gamma)/5.$ Next, suppose that $\delta\leq|\gamma|\leq 2\delta.$ Then
we have $|\frm^{-n_0} \gamma_j|>\frac{1}{2(\lambda+\eps)^{n_0}} |\gamma|$ and the number of curves is at least
$\frac{1}{10}(\lambda-\eps)^{n_0}$ where $\lambda$ is the expansion coefficient of $\hat\frm$ and $\eps$ can be made as
small as we wish by taking $\brI$ large.

Continuing this procedure inductively we construct a Cantor set inside $\Gamma$ such that each interval has 
at least $\frac{1}{10}(\lambda-\eps)^{n_0}$ children and ratios of the lengths of children 
to the length of the parent are at least
$\frac{1}{2(\lambda+\eps)^{n_0}}.$ 
It follows that the resulting Cantor 
set has dimension at least
$$ \frac{\ln \left(\frac{1}{10}(\lambda-\eps)^{n_0}\right)}{\ln \left(2(\lambda+\eps)^{n_0}\right)}. $$
This number can be made as close to $1$ as we wish by taking $n_0$ large and then taking $\brI$
large to make $\eps$ as small as needed.
\end{proof}

\begin{rmk}
The Cantor set above is constructed by taking as children the sub-interval where energy grows by
$\kappa \sqrt{n_0}.$ However, the same estimate remains valid if we take sometimes children with
increasing energy and sometimes children with decreasing energy as long as $I_n$ always stays above
$\brI.$ For example we can require that the energy grows until it reaches $2 \brI$ then decays until
it falls below $\frac{3\brI}{2}$ then grows above $3\brI$ then decays below $\frac{3\brI}{2}$ then
grows above $4\brI$ etc. Then the argument presented above shows that the set of oscillatory orbits
has full Hausdorff dimension. We expect that this set has also positive measure but the proof of this
fact seems out reach at the moment. 
\end{rmk}

\section{Conclusions.}
In this paper we considered piecewise smooth Fermi-Ulam ping pong systems. Near infinity this system can be represented
as a small perturbation of the identity map. Small smooth perturbations of the identity were studied in the context of
inner \cite{L} and outer (\cite{Do}) billiards.  In this case, after a suitable change of coordinates, the problem can
be reduced to the study of small perturbations of the map
\begin{align*}
  \tau_{n+1}&=\tau_n+\omega(\act_n),& \act_{n+1}&=\act_n.
\end{align*}
This map is integrable so the above mentioned problems fall in the context of small smooth
perturbations of integrable systems (i.e. KAM theory). In the case of piecewise smooth perturbations
the normal form also exists: it is a piecewise linear map of a torus. However in contrast with the smooth
case the dynamics of the limiting map is much more complicated and, in fact, it is not completely understood,
especially then the linear part is not hyperbolic. In this paper we described for a simple model example:
\begin{enumerate}[(i)]
\item how to obtain the limiting map and
\item how the properties of the limiting map can be translated to results about
the diffusion for the actual systems.
\end{enumerate}
However, there are plenty of open question on both stages of this procedure.
For example, for piecewise smooth Fermi-Ulam ping pongs it is unknown if there is a positive measure
set of \emph{oscillatory} orbits such that
\begin{align*}
  \lim\inf v_n&<\infty,& \lim\sup v_n&=\infty
\end{align*}
in fact no such orbit is known for $\Delta\in (0,4).$
This demonstrates that more effort is needed in order to develop a general theory of piecewise
smooth near integrable systems.
\begin{figure}[!h]
  \includegraphics[width=8cm]{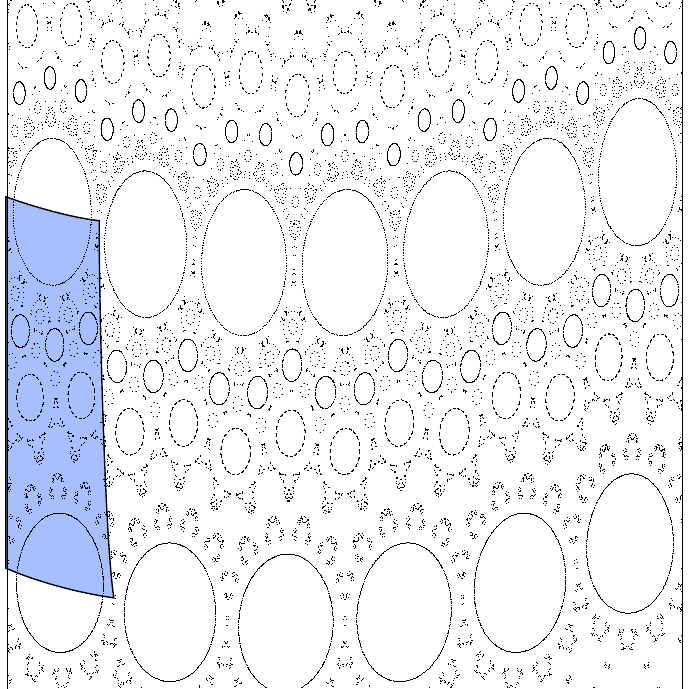}
  \caption{Phase portrait of the region $\torus\times[12,16]$ for selected orbits of the map $f$ where
    $\el(t)=1-0.12\sin(\pi t)$. Notice the similarity of the phase portrait in Figure \ref{f_tfrmElliptic} in the elliptic
    case and the restriction to the shaded area of the phase portrait for the map $f$. The shaded area in fact is a
    fundamental domain of the map $\hat F$}
  \label{f_fElliptic}
\end{figure}
\bibliographystyle{abbrv}
\bibliography{fum}

\end{document}